\documentclass[reqno,11pt]{amsart}

\usepackage{amsmath,amssymb,amsfonts,amsthm, amscd,indentfirst}
\usepackage{amsmath,latexsym,amssymb,amsmath,%awheremscd,
	amscd,amsthm,amsxtra,mathrsfs}
 
\usepackage[top=3cm, bottom=3.5cm, right=2.5cm, left=2.2cm]{geometry}

\usepackage[hyphens]{url} 
\usepackage{hyperref}
\usepackage{bookmark}
\usepackage{amsmath,mathtools}
\usepackage{etoolbox}

\allowdisplaybreaks 
\mathtoolsset{showonlyrefs=true}
\DeclareRobustCommand{\textsection}{\ifmmode\mathsection\else\S\fi}

\usepackage[
  backend=biber,
  style=alphabetic,
  sorting=nyt,
  maxbibnames=99,
  giveninits=true,
]{biblatex}
\renewbibmacro{in:}{}

\DeclareFieldFormat[article,inproceedings,incollection,unpublished]{title}{#1}
\DeclareFieldFormat[article,inproceedings,incollection,unpublished]{citetitle}{#1}

\DeclareSourcemap{%
  \maps[datatype=bibtex]{%
    \map{\step[fieldsource=shortjournal, fieldtarget=journaltitle]}%
  }%
}

\makeatletter
\newcommand*\bbx@lasthash{}
\newtoggle{bbx@dashed}
\AtBeginBibliography{\global\let\bbx@lasthash\empty}
\AtEveryBibitem{%
  \global\togglefalse{bbx@dashed}%
  \iffieldundef{fullhash}
    {\global\let\bbx@lasthash\empty}
    {\iffieldequalstr{fullhash}{\bbx@lasthash}{\global\toggletrue{bbx@dashed}}{}%
     \xdef\bbx@lasthash{\thefield{fullhash}}}%
}
\DeclareNameFormat{labelname}{%
  \ifboolexpr{togl {bbx@dashed}}
    {\ifnum\value{listcount}=1\relax\bibnamedash\fi}
    {\nameparts{#1}%
     \usebibmacro{name:family-given}
       {\namepartfamily}
       {\namepartgiven}
       {\namepartprefix}
       {\namepartsuffix}%
     \usebibmacro{name:andothers}}%
}
\makeatother

\newcounter{mparcnt}

\usepackage{fancyhdr}
\usepackage{esint}
\usepackage{enumerate}
\usepackage{xcolor}

\usepackage{pictexwd,dcpic}
\usepackage{graphicx}
\usepackage{caption}

\usepackage{graphicx}
\usepackage{caption}
\usepackage{slashed}

\usepackage{epsfig,here}
\usepackage{subfigure,here}

\newtheorem{theorem}{Theorem}[section]%[section]
\newtheorem{lemma}[theorem]{Lemma}%[section]
\newtheorem{proposition}[theorem]{Proposition}%[section]
\newtheorem{definition}[theorem]{Definition}%[section]
\newtheorem{remark}[theorem]{Remark}

\newcommand{\abs}[1]{\lvert#1\rvert}
\newcommand{\Abs}[1]{\left\lvert#1\right\rvert}

\newcommand{\norm}[1]{\lVert#1\rVert}
\newcommand{\Norm}[1]{\left\lVert#1\right\rVert}

\newcommand{\rd}{{\rm d}}
\newcommand{\rdV}{{\rm dV}}

\newcommand{\rid}{{\rm id}}
\newcommand{\rtr}{{\rm tr}}
\newcommand{\rdiv}{{\rm div}}

\newcommand{\rRe}{{\rm Re}}

\newcommand{\curl}{{\rm curl}\,}
\newcommand{\ip}{\lrcorner\,}
\newcommand{\w}{\wedge}

\def\<{\langle}
\def\>{\rangle}
\def\S{\mathbb{S}}
\def\R{\mathbb{R}}

\newcommand{\CP}{\mathbb{CP}}

\newcommand{\ra}{\rightarrow}

\newcommand{\eq}[1]{\begin{equation}\allowdisplaybreaks\begin{alignedat}{2} #1 \end{alignedat}\end{equation}}

\numberwithin{equation} {section}

\begin{document}

%------
	
\title[The sharp Chern--Simons--Sobolev inequality and its variants]
{The sharp Chern--Simons--Sobolev inequality and its variants}
\date{\today}

\author{Mingwei Zhang}
\address{ Wuhan University, School of Mathematics and Statistics, 430072 Wuhan, China and 
Albert-Ludwigs-Universit\"at Freiburg,
Mathematisches Institut,
Ernst-Zermelo-Str. 1,
D-79104 Freiburg, Germany}
\email{zhangmwmath@whu.edu.cn}

%------
\begin{abstract}
We introduce and prove a new conformally invariant Chern--Simons--Sobolev inequality, and classify all optimizers, which are the Hopf map $\pi:\mathbb{S}^{2n+1}\to\mathbb{CP}^n$ composed with a symplectic diffeomorphism of $\mathbb{CP}^n$ and an orientation-preserving conformal transformation of $\mathbb{S}^{2n+1}$. The proof uses the sharp curl-Sobolev inequality and a rigidity argument for its equality case.

Motivated by the main result, we also introduce the spherical Radu--Tchrakian--Yang functional, and prove the global
minimality of the Hopf map for radii in a neighborhood of
$((n!)^2/n)^{1/(4n-4)}$, and characterize all minimizers.
We also prove local instability for sufficiently small and sufficiently large radii, with explicit lower bounds for the Morse index.


\end{abstract}

\maketitle
\tableofcontents
%=========

\section{Introduction}

Let $\pi:\S^{2n+1}\to\CP^n$ be the Hopf map. By convention, we equip
$\CP^n$ with the Fubini--Study metric of holomorphic sectional
curvature $4$, so that $\pi$ is a Riemannian submersion. Denote by $\omega$ its K\"ahler form and by $\omega_{2n+1}$ the volume of the standard round sphere $\S^{2n+1}$. For a smooth map $u:\S^{2n+1}\to\CP^n$, the $2$-form $u^*\omega$ is exact, hence admits a potential $1$-form $\alpha$ such that $\rd\alpha=u^*\omega$,
and define
\eq{
    {\rm CS}(u)\coloneqq\int_{\S^{2n+1}}\alpha\w(\rd\alpha)^n.
}
This is independent of the choice of $\alpha$, and is a homotopy invariant, see the discussion in Subsection~\ref{sec2.5}. When $n=1$, ${\rm CS}(u)$ coincides, up to a normalization, with the classical Hopf invariant, see Whitehead's integral formula \cite{Whitehead47}.
In particular, $\operatorname{CS}(\pi)=\pi^{n+1}$.

Our first result concerns the K\"ahler $p$-energy with fixed Chern--Simons invariant ${\rm CS}(u)$ that
\eq{
    \mathcal{E}_p(u)\coloneqq\int_{\S^{2n+1}}\abs{u^*(\omega^p)}^{\frac{2n+1}{2p}}\,\rdV,\qquad 1\leq p\leq n,
}
where $p$ is an integer and $\omega^p$ is the $p$-fold wedge product of $\omega$. The exponent ensures that $\mathcal{E}_p$ is conformally invariant on the domain. Moreover, post-composition by a symplectic diffeomorphism of $\mathbb{CP}^n$ preserves $\mathcal{E}_p$. In particular,
\eq{
    \mathcal{E}_p(\pi)=\left(p!\sqrt{\binom np}\right)^{\frac{2n+1}{2p}}\omega_{2n+1}.
}
It is therefore natural to ask whether the Hopf map $\pi$ minimizes $\mathcal{E}_p$ in its homotopy class, and whether these transformations give all minimizers.

We prove that this is true when $n$ is odd and
$1\leq p\leq \frac{n+1}{2}$.

\begin{theorem}\label{thm1.1}
Let $n\geq1$ be odd and $1\leq p\leq\frac{n+1}{2}$. Then the minimizers of $\mathcal{E}_p(u)$ in the homotopy class of the Hopf map $\pi$ are exactly $h\circ\pi\circ\Phi$, where $h:\CP^n\to\CP^n$ is a symplectic diffeomorphism (i.e. $h^*\omega=\omega$), and $\Phi\in{\rm Conf}^+(\S^{2n+1})$ is an orientation-preserving conformal transformation. Equivalently, we have the sharp Chern--Simons--Sobolev inequality
\eq{
    \frac{ \int_{\S^{2n+1}}\abs{u^*(\omega^p)}^{\frac{2n+1}{2p}}\,\rdV }{ \big( \int_{\S^{2n+1}}\alpha\w(\rd\alpha)^n \big)^{\frac{2n+1}{2n+2}} } \geq \frac{2\sqrt{\pi}}{n!}\left(p!\sqrt{\binom np}\right)^{\frac{2n+1}{2p}},
}
with equality case as above.
\end{theorem}

The proof uses the sharp curl--Sobolev inequality of Wang--Z. \cite{WZ26curl}. For odd $n$, we apply this inequality to the $n$-form
\eq{
    \beta\coloneqq\alpha\w(\rd\alpha)^{\frac{n-1}{2}}.
}
Then $\rd\beta=(\rd\alpha)^{(n+1)/2}$ and $\int\beta\w\rd\beta={\rm CS}(u)$, which gives the required lower bound for $p=(n+1)/2$. The remaining cases follow from Maclaurin's inequality. To characterize equality, we use the Killing form equation to recover the Hopf fibration from $u^*\omega$. The resulting map descends along the Hopf fibers to a symplectic diffeomorphism of $\mathbb{CP}^n$.

The variation role of the Hopf map has been studied for several related models. Rivi\`ere \cite{Riviere98CAG} established the global minimality of the Hopf map in its homotopy class for the $q$-energy $\int_{\S^3}\abs{\rd u}^q$ ($q\geq4$) on $3$-sphere, together with the local minimality for the critical exponent $q=3$.  A different conformally invariant functional was studied by Aratyn--Ferreira--Zimerman \cite{AFZ99}, who constructed explicit solutions with nonzero Hopf invariant. When $n=1$, their energy coincides with our $\mathcal{E}_1$ after stereographic projection. The family $\mathcal{E}_p$ considered here extends this energy to higher-dimensional Hopf fibrations and higher wedge powers of the K\"ahler form. For the higher-dimensional extension of Rivi\`ere's result, we refer to Wang--Z. \cite{WZ26curl}.

Another closely related functional is the symplectic Dirichlet energy $\frac{1}{2}\int \abs{u^*\omega}^2$. Speight--Svensson \cites{SpeightSvensson07,SpeightSvensson11} studied its first and second variations and proved that the Hopf map $\pi:\S^3\to\S^2$ minimizes this energy in its homotopy class. The behavior changes in higher dimensions: Slobodeanu \cites{Slobodeanu08,Slobodeanu15} proved that the Hopf map $\pi:\S^{2n+1}\to\CP^n$ is unstable for this energy when $n\geq2$, and established more general instability results for spherical domains of dimension greater than $4$. This dimensional contrast motivates the study of additional terms that can restore stability and minimality.

Our second result is motivated by the Faddeev--Skyrme type result of Radu--Tchrakian--Yang \cite{RTY13}. For the classical Faddeev--Skyrme model, Vakulenko--Kapitanski \cite{VK79} established a topological lower bound in terms of the Hopf invariant. Lin--Yang \cite{LinYang04} obtained an upper bound with the same growth exponent for the least energy and proved the existence of minimizers in infinitely many homotopy classes on $\R^3$. Radu--Tchrakian--Yang \cite{RTY13} developed higher-dimensional analogues with target $\CP^n$ and a Chern--Simons charge.
For $n\geq2$, they considered the energy
\eq{
    \mathcal{E}_{{\rm RTY}}(u) \coloneqq \int_{\R^{2n+1}}\Big(\abs{u^*\omega}^2+\abs{u^*(\omega^n)}^2\Big),\qquad u:\R^{2n+1}\to\CP^n,
}
and proved a lower bound in terms of a fractional power of the Chern--Simons charge. We consider the corresponding functional on the round sphere of radius $R$. By scaling to the unit sphere, it
can be written as
\eq{
    \mathcal E_{\mathrm{RTY}}^R(u) \coloneqq R^{2n-3}\int_{\S^{2n+1}}\abs{u^*\omega}^2\,\rdV+R^{1-2n}\int_{\S^{2n+1}}\abs{u^*(\omega^n)}^2\,\rdV.
}
The two terms have different scaling powers, and the radius plays a role in the minimizing property of the Hopf map.

\begin{theorem}\label{thm4.1}
Let $n\geq2$. There exists $\epsilon_n\in(0,1/2)$, depending only on $n$, such that the Hopf map $\pi$ minimizes $\mathcal{E}_{{\rm RTY}}^R$ among its homotopy class whenever
\eq{\label{eq1}
    \Big( \frac{(n!)^2}{n}\frac{1-\epsilon_n}{1+\epsilon_n} \Big)^{\frac{1}{4n-4}} \leq R \leq \Big( \frac{(n!)^2}{n}\frac{1+\epsilon_n}{1-\epsilon_n} \Big)^{\frac{1}{4n-4}}.
}
Moreover, assuming \eqref{eq1}, all minimizers are given by $u=h\circ\pi\circ\mathcal{R}$, where $h:\CP^n\to\CP^n$ is a symplectic diffeomorphism, i.e. $h^*\omega=\omega$, and $\mathcal{R}\in{\rm SO}(2n+2)$.
\end{theorem}

\begin{theorem}\label{thm4.2}
Let $n\geq2$. The Hopf map $\pi$ is locally unstable in the following parameter intervals. Precisely,

(1) ${\rm ind}(\pi)\geq2n+2$ if
\eq{\label{eq1.-2}
    R > \Big( \frac{(2n-1)(n!)^2}{2n-3} \Big)^{\frac{1}{4n-4}};
}

(2) ${\rm ind}(\pi)\geq(n-1)(n+1)(n+2)$ if
\eq{\label{eq1.-1}
    0 < R < \Big( \frac{(n!)^2}{4n+3} \Big)^{\frac{1}{4n-4}}.
}
\end{theorem}

To prove the minimizing property, we write
$\mathcal E_{\mathrm{RTY}}^R=R^{2n-3}E_\kappa$, where
$\kappa=R^{4-4n}$ and
\eq{
    E_\kappa(u)\coloneqq\int_{S^{2n+1}}\abs{u^*\omega}^2\,\rdV+\kappa\int_{S^{2n+1}}\abs{(u^*\omega)^n}^2\,\rdV.
}
At $\kappa_*:=n/(n!)^2$, the Chern--Simons constraint and the first eigenvalue estimate for co-closed $1$-forms give a lower bound attained by the Hopf map. We further prove that the energy excess controls the difference between the two terms in $E_{\kappa_*}$. This quantitative estimate allows us to vary $\kappa$ in a neighborhood of $\kappa_*$ and to retain the same equality case. The instability results follow from second variations along conformal flows and along flows generated by suitable horizontal quadratic vector fields.

\medskip

\noindent\textit{Organization of the paper.}
In Section~\ref{sec2}, we recall the properties of differential forms, the curl operator, Killing forms, the Hopf map, and the Chern--Simons invariant used in the proofs. In Section~\ref{sec3}, we prove the minimizing property of the Hopf map for the K\"ahler $p$-energy. Section~\ref{sec4} is devoted to the spherical Radu--Tchrakian--Yang functional and the minimality and instability results.

\section{Preliminaries}\label{sec2}

\subsection{Exterior calculus}

Let $(M,g)$ be an oriented Riemannian manifold of dimension $m\equiv 3\pmod 4$. For each $0\le p\le m$, we denote by $\Omega^p$ the space of smooth $p$-forms and equip it with the $L^2$ inner product induced by $g$. The exterior derivative is denoted by $\rd:\Omega^p\to\Omega^{p+1}$, with the convention $\Omega^{m+1}=0$. We fix the orientation determined by the volume form $\rdV$. The Hodge star operator $*:\Omega^p\to\Omega^{m-p}$ is characterized by
\eq{\label{eq:def_Hodge_star}
    \alpha\w*\beta = \<\alpha, \beta\>\,\rdV.
}
The codifferential, i.e. the $L^2$-adjoint of $\rd$, is
\eq{
    \rd^*\coloneqq (-1)^{m(p+1)+1}*\rd*:\Omega^p\ra\Omega^{p-1}.
}

Since $m\equiv 3\pmod 4$, one has $\rd^*=-*\rd*$ on $\Omega^{\frac{m-1}{2}}$. We recall the standard identities
\eq{
    *^2 = (-1)^{p(m-p)}\,{\rm id},\qquad \rd^2 = (\rd^*)^2=0.
}
In particular, throughout this paper $*^2={\rm id}$.

The Hodge Laplacian is defined by ${\Delta \coloneqq \rd\rd^* + \rd^*\rd}$
and the rough Laplacian ${\nabla^*\nabla \coloneqq -\rtr\nabla^2}$.
These operators are related by the Weitzenb\"ock formula
\eq{
    \Delta = \nabla^*\nabla + \sum_{i,j} e^j \w e_i \ip R(e_i, e_j),
}
where $\{e_i\}$ is a local orthonormal frame, $\{e^i\}$ is the dual $1$-form coframe, $\ip$ is the interior multiplication, and $R$ is the curvature tensor. On the standard round sphere $\S^m$, acting on $p$-forms, this reduces to
\eq{\label{sphere_Weitzenbock}
    \Delta = \nabla^*\nabla + p(m-p).
}

We record several identities that will be used repeatedly.

\begin{proposition}\label{basic_form}
Let $(M,g)$ be an oriented closed Riemannian manifold of dimension $m\equiv 3\pmod 4$. For any $\alpha,\beta\in\Omega^{\frac{m-1}{2}}$ and $X\in\Gamma(TM)$, the following identities hold:
\begin{enumerate}
    \item $\abs{X\ip\alpha}^2+\abs{X^\flat\w\alpha}^2=\abs{X}^2\abs{\alpha}^2$;
    \item $X\ip*\!\alpha=-*(X^\flat\w\alpha)$;
    \item $\<X\ip*\!\alpha,\beta\>=-\<\alpha, X\ip*\beta\>$, and in particular $\<X\ip*\!\alpha,\alpha\>=0$;
    \item $\rd = \sum_i e^i\w\nabla_{e_i}$ and $\rd^*=-\sum_i e_i\ip\nabla_{e_i}$;
    \item $\sum_i e^i\w e_i\ip\alpha = \frac{m-1}{2}\alpha$ and $\sum_i e_i\ip e^i\w\alpha = \frac{m+1}{2}\alpha$;
    \item $e^i\w e_j\ip\alpha + e_j\ip e^i\w\alpha = \delta_{ij}\alpha$.
\end{enumerate}
More generally, for $\alpha\in\Omega^p$ and $\beta\in\Omega^{m-p-1}$ the above formulas remain valid, with the following modifications:
\begin{enumerate}
    \item[(2)] $X\ip*\!\alpha=(-1)^p*(X^\flat\w\alpha)$;
    \item[(3)] $\<X\ip*\!\alpha,\beta\>=(-1)^p\<\alpha, X\ip*\beta\>$;
    \item[(5)] $\sum_i e^i\w e_i\ip\alpha = p\alpha$, \quad $\sum_i e_i\ip e^i\w\alpha = (m-p)\alpha$.
\end{enumerate}
In particular,
\eq{\label{direct_dual}
    \<X\ip\alpha,\beta\> = \<\alpha,X^\flat\w\beta\>, \quad\forall \alpha\in\Omega^p,\ \beta\in\Omega^{p-1}.
}
\end{proposition}

\subsection{The curl operator}\label{sec2.2}

Let $m\equiv 3\pmod 4$ and set $p\coloneqq\frac{m-1}{2}$. We define the curl operator by
\eq{\label{eq:curl-def}
    \curl \coloneqq *\rd: \Omega^{\frac{m-1}{2}}\ra \Omega^{\frac{m-1}{2}}.
}
When $m=3$ and $\alpha$ is the $1$-form dual to a vector field $X$, definition \eqref{eq:curl-def} coincides with the usual vector-calculus curl.

A basic computation using $\<\alpha,\beta\>\rdV=\alpha\wedge*\beta$ and Stokes' theorem shows that on a closed manifold $M$ the curl operator is formally self-adjoint on $\Omega^{\frac{m-1}{2}}(M)$:
\eq{\label{eq:curl-selfadj}
    \int_M \< \curl\alpha,\beta\>\,\rdV
    = \int_M \< \alpha,\curl\beta\>\,\rdV.
}
In particular,
\eq{\label{eq:curl-energy}
    \int_M \< \curl\alpha,\alpha\>\,\rdV = \int_M \alpha\wedge \rd\alpha.
}

Let $\tilde g=\rho^2g$ for some smooth function $\rho>0$. The standard conformal covariance is
\eq{\label{conformal_change}
    *_{\tilde{g}} = \rho^{m-2p} *_{g}, \quad \<\cdot,\cdot\>_{\tilde{g}} = \rho^{-2p}\<\cdot,\cdot\>_g \quad\hbox{on}\ \Omega^p.
}
In particular, for any $\alpha\in\Omega^{\frac{m-1}{2}}$,
\eq{
    \curl_{\tilde{g}}\alpha = *_{\tilde{g}}\rd\alpha = \rho^{-1}\,\curl_g\alpha.
}
Together with \eqref{conformal_change} this yields
\eq{
    \abs{\curl_{\tilde{g}}\alpha}_{\tilde{g}}^{\frac{2m}{m+1}} = \rho^{-m}\abs{\curl_g\alpha}_g^{\frac{2m}{m+1}},\qquad
    \<\curl_{\tilde{g}}\alpha, \alpha\>_{\tilde{g}} = \rho^{-m}\<\curl_g\alpha,\alpha\>_g.
}

\subsection{Killing forms}
We recall the definition and basic properties of Killing forms.

\begin{definition}\label{def_Killing}
A co-closed $p$-form $\alpha\in\Omega^p$ is called a \emph{Killing $p$-form} if
\eq{\label{eq:104}
    \nabla_X\alpha = \frac{1}{p+1}X\ip\rd\alpha,\quad \forall X\in\Gamma(TM).
}
\end{definition}

In particular, an $\frac{m-1}{2}$-Killing form is defined by
\eq{
    \nabla_X\alpha = \frac{2}{m+1}X\ip\rd\alpha,\quad \forall X\in\Gamma(TM).
}
On $\S^m$, each $\frac{m-1}{2}$-Killing form $\alpha$ can be decomposed as
\eq{
    \alpha = \alpha_++\alpha_-,
}
where $\alpha_\pm$ is a first $\pm$-eigenform of $\curl$, i.e.
\eq{
    \curl\alpha_+=\frac{m+1}{2}\alpha_+, \quad \curl\alpha_-=-\frac{m+1}{2}\alpha_-.
}
We call $\alpha$ a positive (respectively negative) $\frac{m-1}{2}$-Killing form if $\alpha=\alpha_+$ (respectively $\alpha=\alpha_-$).

The following properties are important; analogous statements hold for negative $\frac{m-1}{2}$-Killing forms.

\begin{proposition}\label{Killing_properties}
Let $\alpha$ be a positive $\frac{m-1}{2}$-Killing form on $\S^m$. Then
\begin{enumerate}
    \item $\nabla_X\alpha = X\ip*\alpha$ for every vector field $X$;
    \item $\abs{\alpha}$ is constant.
\end{enumerate}
\end{proposition}

\begin{proof}
(1) By definition we have $\curl\alpha = \frac{m+1}{2}\alpha$. Combining with Definition~\ref{def_Killing}, we have
\eq{
    \nabla_X\alpha = \frac{2}{m+1}X\ip \rd\alpha = \frac{2}{m+1}X\ip *\curl\alpha = X\ip *\alpha.
}

(2) It follows from (1) and Proposition~\ref{basic_form}(3).
\end{proof}

\subsection{The Hopf map}

We view the sphere as
\eq{
    \S^{2n+1} = \{z=(z_0,z_1,\cdots,z_n)\in\mathbb{C}^{n+1}\,:\, \abs{z}=1\}.
}
We agree with the convention that the complex projective space $\mathbb{CP}^n$ is equipped with the Fubini--Study metric
\eq{
    g_{i\bar{j}} = \frac{(1+\abs{z}^2)\delta_{i\bar{j}} - \bar{z}_iz_j}{(1+\abs{z}^2)^2}.
}
The standard \textit{Hopf map} (or \textit{Hopf fibration}) $\pi:\S^{2n+1}\to\mathbb{CP}^n$ is defined by
\eq{
    \pi(z) \coloneqq [z_0:z_1:\cdots:z_n].
}
The choice of the Fubini--Study metric makes the Hopf map a smooth Riemannian submersion, and also a harmonic map.

By definition, the Hopf map $\pi$ projects every complex line $\mathbb{C}z$ to $[z]$. Therefore it can be viewed as the quotient map of the principal $\S^1$-bundle
\eq{
    \S^1 \rightarrow \S^{2n+1} \xrightarrow{\pi} \mathbb{CP}^n,
}
where $\S^1$ freely acts on $\S^{2n+1}$ by
\eq{
    e^{i\theta}\cdot z \coloneqq e^{i\theta}z.
}
Each $\S^1$-fiber is an orbit of the flow induced by the unit Reeb field $\xi\in\ker(\rd\pi)$.

\subsection{The Chern--Simons invariant}\label{sec2.5}

Let $u:\S^{2n+1}\to\CP^n$ be a smooth map. We equip $\CP^n$ with the Fubini--Study metric of holomorphic sectional curvature $4$, so that the Hopf map $\pi:\S^{2n+1}\to\CP^n$ is a Riemannian submersion. Denote by $\omega$ the K\"ahler form of $\CP^n$, and
\eq{
    \omega^p \coloneqq \omega\w\cdots\w\omega \quad\hbox{($p$ times)}.
}
Since the K\"ahler form $\omega$ is closed, so is $u^*\omega\in\Omega^2(\S^{2n+1})$. As $H^2(\S^{2n+1})=0$, it is exact, so we may choose a $1$-form $\alpha\in\Omega^1(\S^{2n+1})$ such that
\eq{
    u^*\omega = \rd\alpha.
}
By the Hodge theory, $\alpha$ can be chosen to be co-closed.

We consider the Chern--Simons type invariant
\eq{
    {\rm CS}(u) \coloneqq \int_{\S^{2n+1}} \alpha\wedge(\rd\alpha)^n, \quad \alpha\in\Omega^1(\S^{2n+1}).
}
It is independent of the choice of the potential $\alpha$, since
\eq{
    \int_{\S^{2n+1}} \alpha\wedge(\rd\alpha)^n = \int_{\S^{2n+1}} (\alpha+\rd\phi)\wedge(\rd(\alpha+\rd\phi))^n
}
for any function $\phi$.

We now show that ${\rm CS}(u)$ is a homotopy invariant.

\begin{lemma}
For any $u:\S^{2n+1}\to\mathbb{CP}^n$, there exists a lift-up $\tilde{u}:\S^{2n+1}\to\S^{2n+1}$ such that
\eq{
    u = \pi\circ\tilde{u}.
}
\end{lemma}
\begin{proof}
Given $x\in\S^{2n+1}$, $u(x)\in\mathbb{CP}^n$ is viewed as a complex line in $\mathbb{C}^{n+1}$, hence $u$ yields a complex line bundle $L_u$ on $\S^{2n+1}$. 

We claim that one may choose $\tilde{u}(x)$ to be a unit vector on the line $u(x)$, hence the conclusion follows. It suffices to show that $L_u$ admits a nowhere vanishing section. In fact, the line bundles are characterized by the first Chern class, while
\eq{
    c_1(L_u) \in H^2(\S^{2n+1},\mathbb{Z}) = 0,
}
hence $L_u$ is trivial. The claim follows.
\end{proof}

In particular, for $u=\pi$, let $\eta\in\Omega^1(\S^{2n+1})$ be defined by
\eq{
    \eta=\frac{1}{2}\sum_j(x_j\rd y_j-y_j\rd x_j),\qquad \pi^*\omega = \rd\eta,
}
where $\pi$ is the Hopf map. In fact, $\theta\coloneqq 2\eta$ is the dual $1$-form of the Reeb field. Let $\{e_1,\cdots,e_{2n},\theta\}$ be an orthonormal basis, then
\eq{\label{eq:1}
    \rd\eta = \pi^*\omega = \sum_{j=1}^n e^{2j-1}\w e^{2j}.
}
It follows that
\eq{
    (\rd\eta)^n = n!e^1\w\cdots\w e^{2n},
}
and hence
\eq{
    \eta\w(\rd\eta)^n = \frac{n!}{2}\rdV.
}
Therefore
\eq{\label{eq0}
    {\rm CS}(\pi) = \int_{\S^{2n+1}} \eta\wedge(\rd\eta)^n = \frac{n!}{2}\omega_{2n+1} = \pi^{n+1}.
}

In order to obtain a integer-valued invariant, we define the normalized invariant
\eq{
    Q(u) \coloneqq \frac{{\rm CS}(u)}{{\rm CS}(\pi)} = \frac{1}{\pi^{n+1}}\int_{\S^{2n+1}} \alpha\wedge(\rd\alpha)^n, \quad \alpha\in\Omega^1(\S^{2n+1}).
}

\begin{proposition}\label{prop9.3}
${\rm CS}(u)$ is a homotopy invariant, i.e. ${\rm CS}(u)={\rm CS}(v)$ iff $u$ is homotopic to $v$ (denoted by $u\simeq v$). Moreover,
\eq{
    {\rm CS}(u) = \deg(\tilde{u}){\rm CS}(\pi) = \pi^{n+1}\deg(\tilde{u}).
}
Consequently,
\eq{
    Q: [\S^{2n+1},\mathbb{CP}^n]=\pi_{2n+1}(\mathbb{CP}^n) \to \mathbb{Z}, \qquad Q(u)=\deg(\tilde{u})
}
is bijective.
\end{proposition}
\begin{proof}
Since $u=\pi\circ\tilde{u}$, we have
\eq{
    u^*\omega = \tilde{u}^*\pi^*\omega = \tilde{u}^*\rd\eta = \rd(\tilde{u}^*\eta).
}
Hence $\alpha=\tilde{u}^*\eta$ is a canonical choice of the potential $1$-form. It follows
\eq{
    {\rm CS}(u) = \int_{\S^{2n+1}} \tilde{u}^*(\eta\w(\rd\eta)^n) = \deg(\tilde{u})\int_{\S^{2n+1}} \eta\w(\rd\eta)^n = \deg(\tilde{u}){\rm CS}(\pi).
}
Therefore
\eq{
    {\rm CS}(u)={\rm CS}(v) \iff \tilde{u} \simeq \tilde{v}.
}
It remains to show that
\eq{
    \tilde{u} \simeq \tilde{v} \iff u \simeq v,
}
or equivalently, the induced map $\pi_*$ gives a bijection
\eq{
    \pi_*:[\S^{2n+1},\S^{2n+1}] \to [\S^{2n+1},\mathbb{CP}^n].
}
In fact, if $\tilde{u}\simeq\tilde{v}$, then obviously $u\simeq v$. Reversely, the long exact sequence
\eq{
    0 = \pi_{2n+1}(\S^1) \to \pi_{2n+1}(\S^{2n+1}) \overset{\pi_*}{\to} \pi_{2n+1}(\CP^n) \to \pi_{2n}(\S^1) = 0
}
implies that $\pi_*$ is an isomorphism.
\end{proof}

\begin{remark}
In the case $n=1$, $Q(u)$ is the classical Hopf invariant.
\end{remark}

\section{The Hopf map minimizes the K\"ahler \texorpdfstring{$p$}{p}-energy}\label{sec3}

Define the K\"ahler p-energy of $u$ by
\eq{
    \mathcal{E}_p(u) \coloneqq \int_{\S^{2n+1}} \abs{u^*(\omega^p)}^\frac{2n+1}{2p}.
}
By \eqref{conformal_change}, it is conformally invariant. In particular, for $u=\pi$, by \eqref{eq:1} we have
\eq{
    \pi^*(\omega^p) = p!\sum_{1\leq j_1<\cdots<j_p\leq n}(e^{2j_1-1}\w e^{2j_1})\w\cdots\w(e^{2j_p-1}\w e^{2j_p}),
}
and hence
\eq{
    \abs{\pi^*(\omega^p)}^2 = (p!)^2\binom{n}{p}.
}
Therefore
\eq{
    \mathcal{E}_p(\pi) = \Bigg( p! \sqrt{\binom{n}{p}} \Bigg)^{\frac{2n+1}{2p}}\omega_{2n+1}.
}

We now prove Theorem~\ref{thm1.1}.

\begin{proof}[Proof of Theorem~\ref{thm1.1}]
The $n=1$ case was proved in \cite{WZ26curl}. Here we handle $n>1$ cases. 
We first prove the case $p=\frac{n+1}{2}$. Set
\eq{\label{eq:96}
    \beta \coloneqq \alpha\w(\rd\alpha)^{\frac{n-1}{2}} \in\Omega^n(\S^{2n+1}).
}
Then
\eq{\label{eq:97}
    \rd\beta = (\rd\alpha)^{\frac{n+1}{2}} \in\Omega^{n+1}(\S^{2n+1}).
}
\cite[Theorem~1.1]{WZ26curl} implies that
\eq{\label{eq:98}
    \Big( \int_{\S^{2n+1}} \abs{\rd\beta}^{\frac{2n+1}{n+1}} \Big)^{\frac{2n+2}{2n+1}} \geq (n+1)\omega_{2n+1}^{\frac{1}{2n+1}}\int_{\S^{2n+1}}\beta\w\rd\beta,
}
with equality if and only if $\beta$ is a positive Killing $n$-form modulo orientation-preserving conformal transformations and modulo $\ker(\rd)$.

On one hand, by \eqref{eq:96} and \eqref{eq:97} we have
\eq{\label{eq:99}
    \int_{\S^{2n+1}}\beta\w\rd\beta = \int_{\S^{2n+1}}\alpha\w(\rd\alpha)^n = {\rm CS}(u).
}
Therefore, by Proposition~\ref{prop9.3}, $\int\beta\w\rd\beta={\rm CS}(u)$ is fixed in the homotopy class of the Hopf map $\pi$.

On the other hand, since $\rd\alpha=u^*\omega$, we have
\eq{
    \rd\beta = (\rd\alpha)^{\frac{n+1}{2}} = u^*(\omega^{\frac{n+1}{2}}).
}
Hence
\eq{\label{eq:100}
    \int_{\S^{2n+1}} \abs{\rd\beta}^{\frac{2n+1}{n+1}} = \mathcal{E}_{\frac{n+1}{2}}(u).
}
Combining \eqref{eq:98}, \eqref{eq:99} and \eqref{eq:100} yields
\eq{\label{eq:103}
    \mathcal{E}_{\frac{n+1}{2}}(u) \geq \Big( (n+1)\omega_{2n+1}^{\frac{1}{2n+1}} {\rm CS}(u) \Big)^{\frac{2n+1}{2n+2}}.
}

We now consider the general case $1\leq p\leq\frac{n+1}{2}$. At a given point $x\in\S^{2n+1}$, we choose an orthonormal basis $\{e_j\}_{j=1}^{2n+1}$ such that
\eq{
    \rd\alpha(x) = u^*\omega(x) = \sum_{j=1}^n\lambda_j(x)e^{2j-1}(x)\w e^{2j}(x), \quad \lambda_j=\lambda_j(x)\geq0.
}
It follows
\eq{\label{eq:101}
    \abs{(\rd\alpha)^p}^2 = (p!)^2\sigma_p(\lambda_1^2,\cdots,\lambda_n^2),
}
where $\sigma_p$ is the $p$-th elementary symmetric polynomial. The Maclaurin's inequality implies that
\eq{\label{eq:102}
    \left( \frac{\sigma_p(\lambda_1^2,\cdots,\lambda_n^2)}{ \binom{n}{p} } \right)^{\frac{1}{p}} \geq \left( \frac{\sigma_{\frac{n+1}{2}}(\lambda_1^2,\cdots,\lambda_n^2)}{ \big(\!\substack{n\\[1pt]\frac{n+1}{2}}\!\big) } \right)^{\frac{2}{n+1}}.
}
When $1\leq p<\frac{n+1}{2}$, equality holds if and only if $\lambda_1=\cdots=\lambda_n$ or the number of non-zero $\lambda_j$'s is smaller than $p$. Combining \eqref{eq:101} and \eqref{eq:102} gives
\eq{
    \left( \frac{\abs{(\rd\alpha)^{p}}}{p!\sqrt{ \binom{n}{p} }} \right)^{\frac{1}{p}} \geq \left( \frac{\abs{(\rd\alpha)^{\frac{n+1}{2}}}}{(\frac{n+1}{2})!\sqrt{ \big(\!\substack{n\\[1pt]\frac{n+1}{2}}\!\big) }} \right)^{\frac{2}{n+1}}
}
Hence
\eq{
     \abs{(\rd\alpha)^{p}}^{\frac{2n+1}{2p}} \geq \left(p!\sqrt{\binom{n}{p}}\,\right)^{\frac{2n+1}{2p}} \left( \frac{\abs{(\rd\alpha)^{\frac{n+1}{2}}}}{(\frac{n+1}{2})!\sqrt{ \big(\!\substack{n\\[1pt]\frac{n+1}{2}}\!\big) }} \right)^{\frac{2n+1}{n+1}}.
}
Integrating and using \eqref{eq:103} give
\eq{
    \mathcal{E}_p(u) \geq C(n,p) \mathcal{E}_{\frac{n+1}{2}}(u) \geq C(n,p)\Big( (n+1)\omega_{2n+1}^{\frac{1}{2n+1}} {\rm CS}(u) \Big)^{\frac{2n+1}{2n+2}}.
}
where
\eq{
    C(n,p) = \frac{ \left(p!\sqrt{\binom{n}{p}}\,\right)^{\frac{2n+1}{2p}} }{ \left( (\frac{n+1}{2})!\sqrt{ \big(\!\substack{n\\[1pt]\frac{n+1}{2}}\!\big) } \right)^{\frac{2n+1}{n+1}} }.
}

Moreover, equality holds if and only if 
\begin{enumerate}
    \item equality in \eqref{eq:98} holds, i.e. $\beta$ is a positive Killing $n$-form modulo orientation-preserving conformal transformations and modulo $\ker(\rd)$.;
    \item equality in \eqref{eq:102} holds.
\end{enumerate}

We now claim that $\lambda(x)\equiv\lambda$ is constant for every $x\in\S^{2n+1}$.

\begin{proof}[Proof of Claim]
After modulo an orientation-preserving conformal transformation $\Phi\in{\rm Conf}^+(\S^{2n+1})$, (1) implies that $\beta = \beta_0 + \gamma$ for some $\gamma\in\ker(\rd)$, where $\beta_0$ is a positive Killing $n$-form with $*\rd\beta_0=(n+1)\beta_0$. Hence $\rd\beta=\rd\beta_0=(\rd\alpha)^{\frac{n+1}{2}}$. It follows
\eq{
    *(\rd\alpha)^{\frac{n+1}{2}} = *\rd\beta = *\rd\beta_0 = (n+1)\beta_0.
}
Differentiating gives
\eq{\label{eq:105}
    \nabla_X\beta_0 = \frac{1}{n+1}*\nabla_X((\rd\alpha)^{\frac{n+1}{2}}) = \frac{1}{2}*(\nabla_X(\rd\alpha)\w(\rd\alpha)^{\frac{n-1}{2}})
}
for any vector field $X$. Since $\beta_0$ is a positive Killing $n$-form, by definition \eqref{eq:104} we have
\eq{\label{eq:106}
    \nabla_X\beta_0 = \frac{1}{n+1}X\ip\rd\beta_0 = \frac{1}{n+1}X\ip(\rd\alpha)^{\frac{n+1}{2}} = \frac{1}{2}(X\ip\rd\alpha)\w(\rd\alpha)^{\frac{n-1}{2}}.
}
Combining \eqref{eq:105} and \eqref{eq:106} yields
\eq{\label{eq:107}
    *(\nabla_X(\rd\alpha)\w(\rd\alpha)^{\frac{n-1}{2}}) = (X\ip\rd\alpha)\w(\rd\alpha)^{\frac{n-1}{2}}.
}

Set
\eq{\label{eq:113}
    U \coloneqq \{ x\in\S^{2n+1} \mid (\rd\alpha)^n(x) \not=0 \}.
}
Clearly $U\not=\emptyset$. On $U$, $\rd\alpha$ has the maximal rank $2n$, hence $\ker(\rd\alpha)$ has rank $1$. Choose a unit vector field $\xi\in\ker(\rd\alpha)$, and denote $H\coloneqq \xi^\perp$. We now work on $U$.

Given a point $x\in\S^{2n+1}$, for an orthonormal basis $\{e_1,\cdots,e_{2n},\xi\}$, we may write the standard form
\eq{\label{eq:112}
    \rd\alpha(x) = \sum_{j=1}^n \lambda_j(x)e^{2j-1}(x)\w e^{2j}(x), \qquad \lambda_j(x)>0.
}
Decompose
\eq{
    \nabla_X(\rd\alpha) = \nabla_X^H(\rd\alpha) + \xi^\flat\w a_X,
}
where
\eq{
    \nabla_X^H(\rd\alpha)\in \Lambda^2H^*, \qquad a_X\in H^*.
}

(i) Choose $X\in H$. Comparing the $(2n+1)$-th component in both sides of \eqref{eq:107}, we see that
\eq{
    \nabla_X^H(\rd\alpha)\w(\rd\alpha)^{\frac{n-1}{2}} = 0.
}
Since $n-2\geq\frac{n-1}{2}$, it follows
\eq{\label{eq:108}
    \nabla_X^H(\rd\alpha)\w(\rd\alpha)^{n-2} = 0.
}
Since $\rd\alpha|_H=u^*\omega|_H$ is a symplectic form on the $2n$-dimensional space $H$, we may apply the hard Lefschetz theorem (see for instance \cite[Section 6]{Voisin02Book}), which says
\eq{
    L^{n-2}:\Lambda^2H^* \to \Lambda^{2n-2}H^*, \quad \phi\mapsto\phi\w(\rd\alpha)^{n-2}
}
is an isomorphism. Therefore, \eqref{eq:108} implies that
\eq{
    \nabla_X^H(\rd\alpha) = 0, \qquad\forall X\in H,
}
hence
\eq{\label{eq:109}
    \nabla_X(\rd\alpha) = \xi^\flat\w a_X \quad\hbox{with}\quad a_X\in H^*, \qquad\forall X\in H.
}

(ii) Fix $j$ and choose $X=e_{2j-1}$. Then \eqref{eq:107} and \eqref{eq:109} imply that
\eq{\label{eq:114}
    *_H(a_{e_{2j-1}}\w(\rd\alpha)^{\frac{n-1}{2}}) = \lambda_j e^{2j}\w(\rd\alpha)^{\frac{n-1}{2}}.
}
Comparing the $(2j-1,2j)$-component of both sides, we see that
\eq{\label{eq:110}
    a_{e_{2j-1}} = c_j e^{2j-1}
}
for some $c_j$. Set $I_j\coloneqq\{1,\cdots,n\}\backslash\{j\}$. For $K\subset I_j$, denote
\eq{
    \alpha_K \coloneqq \bigwedge_{k\in K}(e^{2k-1}\w e^{2k}), \qquad \lambda_K \coloneq \prod_{k\in K}\lambda_k.
}
Then the right-hand side of \eqref{eq:114} is
\eq{
    \lambda_j e^{2j}\w(\rd\alpha)^{\frac{n-1}{2}} = \lambda_j\cdot\frac{n-1}{2}!\cdot\sum_{\substack{K\subset I_j\\\abs{K}=\frac{n-1}{2}}} \lambda_K e^{2j}\w \alpha_K
}
The left-hand side of \eqref{eq:114} is
\eq{
    *_H(a_{e_{2j-1}}\w(\rd\alpha)^{\frac{n-1}{2}}) = c_j\cdot\frac{n-1}{2}!\cdot\sum_{\substack{K\subset I_j\\\abs{K}=\frac{n-1}{2}}} \lambda_{I_j\backslash K} e^{2j}\w \alpha_K
}
Comparing the coefficient, we obtain that
\eq{\label{eq:115}
    c_j\lambda_{I_j\backslash K} = \lambda_j \lambda_K, \qquad\hbox{for any}\quad K\subset I_j,\quad \abs{K}=\frac{n-1}{2}.
}
Given any $j_1,j_2$, choose $K\subset I_j$, $\abs{K}=\frac{n-1}{2}$ such that $j_1\in K$, $j_2\not\in K$. Set $K'\coloneqq(K\backslash\{j_1\})\cup\{j_2\}$. Then using \eqref{eq:115} we have
\eq{
    c_j\lambda_{I_j\backslash K} = \lambda_j \lambda_K, \qquad c_j\lambda_{I_j\backslash K'} = \lambda_j \lambda_{K'}.
}
It follows
\eq{
    \frac{\lambda_{j_1}}{\lambda_{j_2}} = \frac{\lambda_j \lambda_K}{\lambda_j \lambda_{K'}} = \frac{c_j\lambda_{I_j\backslash K}}{c_j\lambda_{I_j\backslash K'}} = \frac{\lambda_{j_2}}{\lambda_{j_1}}.
}
Hence $\lambda_{j_1}=\lambda_{j_2}$. Therefore, we have
\eq{
    \lambda_1=\cdots=\lambda_n\eqcolon\lambda.
}
Then \eqref{eq:115} implies that
\eq{
    c_1=\cdots=c_n=\lambda.
}
Similarly, for $X=e_{2j}$ we have
\eq{\label{eq:111}
    a_{e_{2j}} = \lambda_j e^{2j} = \lambda e^{2j}.
}
Combining \eqref{eq:109}, \eqref{eq:110}, and \eqref{eq:111} gives
\eq{\label{eq111.1}
    \nabla_X(\rd\alpha) = \lambda \xi^\flat\w X^\flat, \quad\forall X\in H \qquad\hbox{at a point }x.
}

(iii) Choose $X=\xi$. Then the right-hand side of \eqref{eq:107} vanishes, hence
\eq{
    \nabla_{\xi}(\rd\alpha)\w(\rd\alpha)^{\frac{n-1}{2}} = 0.
}
Apply again the hard Lefschetz theorem for the $\Lambda^2H^*$ part and the $\xi^\flat\w H^*$ part respectively, we obtain
\eq{
    \nabla_{\xi}(\rd\alpha) = 0 \qquad\hbox{at a point }x.
}

We now claim that $\lambda(x)$ is constant. In fact, we have for $X\in H$
\eq{
    X(\abs{\rd\alpha}^2) = 2\<\nabla_X(\rd\alpha),\rd\alpha\> = 2\lambda\<\xi^\flat\w X^\flat,\rd\alpha\> = 0,
}
and for $X=\xi$,
\eq{
    \xi(\abs{\rd\alpha}^2) = 2\<\nabla_{\xi}(\rd\alpha),\rd\alpha\> = 0.
}
Since $\abs{\rd\alpha}^2 = n\lambda^2$, we see that $\rd\lambda=0$, hence $\lambda(x)\equiv\lambda$ is constant on the maximal connected component of $U$. Thus the maximal connected component is a closed set. But by definition \eqref{eq:113}, it is an open set. Hence $U=\S^{2n+1}$ since $\S^{2n+1}$ is connected.
\end{proof}

Therefore, $\xi$ is globally defined. Let $g(AX,Y)=\lambda^{-1}\rd\alpha(X,Y)$ with $A$ anti-symmetric. Differentiating $\rd\alpha(\xi,Y)=0$ gives
\eq{
    \nabla_X\xi=AX, \qquad \rd\xi^\flat = \frac{2}{\lambda}\rd\alpha.
}
Hence $\xi$ is a unit Killing vector field. Consequently, $\lambda=1$, and $\lambda^{-1}\rd\alpha$ is a Sasakian $2$-form. Therefore, $u$ descends to a symplectic diffeomorphism, and we complete the proof. For details, see the proof of equality case of \cite[Theorem~1.4]{WZ26curl}.
\end{proof}

\begin{remark}
In particular, when $n=1$, it covers the case $\rho=0$ in the Faddeev--Skyrme model, see \cite[Section~7]{WZ26curl}. In view of AM--GM, it also covers the $3$-energy model, see \cite[Section~6]{WZ26curl}.
\end{remark}

\begin{remark}
We remark that Radu--Tchrakian--Yang \cite{RTY13} considered a similar Faddeev--Skyrme type problem: to minimize
\eq{
    \mathcal{E}_{{\rm RTY}}(u) = \int_{\R^{2n+1}} \Big( \abs{u^*\omega}^2 + \abs{u^*(\omega^n)}^2 \Big), \qquad u:\R^{2n+1}\to\CP^n, \qquad n>1.
}
They proved that there exists some $c>0$ such that
\eq{
    \mathcal{E}_{{\rm RTY}}(u) \geq c\,\abs{{\rm CS}(u)}^{\frac{2n+1}{2n+2}}.
}
\end{remark}

\section{The spherical Radu--Tchrakian--Yang functional}\label{sec4}

Spherical versions of Skyrme-type models provide a natural setting in which the size of the domain becomes a variational parameter. Manton--Ruback \cite{MantonRuback86} studied the Skyrme model on curved spatial domains, including the round $3$-sphere. For the Faddeev–-Skyrme model on $\S^3$, Ward \cite{Ward99} examined the Hopf map and showed that it becomes unstable when the radius is sufficiently large. More recently, Guerra--Lamy--Zemas \cite{GLZ26} proved global minimality of the Hopf map in its homotopy class for sufficiently large coupling. Motivated by these developments, we consider the Radu-–Tchrakian–-Yang functional for maps from the round sphere $\S^{2n+1}_R$ of radius $R$ into $\CP^n$ with $n\geq2$, and investigate how its minimizing properties depend on $R$.

We consider
\eq{
    \mathcal{E}_{{\rm RTY}}(u) = \int_{\S^{2n+1}_R} \Big( \abs{u^*\omega}^2 + \abs{u^*(\omega^n)}^2 \Big), \qquad u:\S^{2n+1}_R\to\CP^n.
}
By scaling we consider equivalently
\eq{
    \mathcal{E}_{{\rm RTY}}^R(u) = R^{2n-3}\int_{\S^{2n+1}} \abs{u^*\omega}^2 + R^{1-2n}\int_{\S^{2n+1}}\abs{u^*(\omega^n)}^2, \qquad u:\S^{2n+1}\to\CP^n.
}

\subsection{Proof of Theorem~\ref{thm4.1}}

As mentioned above, we choose the potential $1$-form $\alpha$ to be co-closed.
For simplicity, we denote $F\coloneqq \rd\alpha=u^*\omega$, and set
\eq{
    E_\kappa(u) \coloneqq \int_{\S^{2n+1}} \abs{F}^2 + \kappa\int_{\S^{2n+1}}\abs{F^n}^2, \qquad \kappa>0.
}
It is clear that
\eq{\label{eq1.0}
    \mathcal{E}_{{\rm RTY}}^R(u) = R^{2n-3}E_{R^{4-4n}}(u).
}
Since
\eq{
    \abs{\pi^*\omega}^2=n, \qquad \abs{(\pi^*\omega)^n}=n!,
}
we see that
\eq{\label{eq1.0.1}
    E_\kappa(\pi) = \big(n+\kappa(n!)^2\big)\omega_{2n+1}.
}

\begin{proposition}
Let $u\cong\pi$ and $\kappa^* \coloneqq n/(n!)^2$. Then there exists an $\epsilon_n\in(0,1/2)$, depending only on $n$, such that 
\eq{\label{eq1.1}
    E_{\kappa_*}(u) - 2n\omega_{2n+1} \geq 2\epsilon_n\Abs{ \int \abs{F}^2 - \kappa_*\int \abs{F^n}^2 }.
}
\end{proposition}
\begin{proof}
First of all, it suffices to consider the case that
\eq{\label{eq1.2}
    E_{\kappa_*}(u) - 2n\omega_{2n+1}\leq1.
}
In fact, if otherwise, then
\eq{
    \Abs{ \int \abs{F}^2 - \kappa_*\int \abs{F^n}^2 } \leq \int \abs{F}^2 + \kappa_*\int \abs{F^n}^2 = E_{\kappa_*}(u) \leq (1+2n\omega_{2n+1})(E_{\kappa_*}(u) - 2n\omega_{2n+1}),
}
as claimed. Hence we may assume \eqref{eq1.2}.

Using \eqref{eq0} we have
\eq{
    \int \Big\<\alpha,\frac{*F^n}{n!}\Big\> = \frac{1}{n!}\int \alpha\w F^n = \frac{1}{2}\omega_{2n+1}.
}
Therefore,
\eq{
    n\Norm{\frac{*F^n}{n!}-2\alpha}_2^2 = \frac{n}{(n!)^2}\norm{F^n}_2^2 + 4n\norm{\alpha}_2^2 - 4n\int\Big\< \alpha,\frac{*F^n}{n!} \Big\> = \kappa_* \norm{F^n}_2^2 + 4n\norm{\alpha}_2^2 - 2n\omega_{2n+1}.
}
It follows that
\eq{\label{eq2}
    E_{\kappa_*}(u) - 2n\omega_{2n+1} = \norm{F}_2^2 - 4n\norm{\alpha}_2^2 + n\Norm{\frac{*F^n}{n!}-2\alpha}_2^2.
}

For later use, we recall a spectrum lemma, see Ikeda--Taniguchi \cite{IT78}: the eigenvalue of the Hodge Laplacian for co-closed $1$-forms on $\S^{2n+1}$ starts as $4n, 6n+3, \cdots$, and the first eigenspace is formed by all Killing $1$-forms.

We now decompose $\alpha=\xi+\beta$, where $\xi$ is a Killing $1$-form and $\beta$ is $L^2$-orthogonal to the space of Killing $1$-forms. Note that $\xi$, as a Killing $1$-form, is a first eigenform of the Hodge Laplacian on the space of co-closed forms, with respect to the eigenvalue $4n$. That is,
\eq{\label{eq2.1}
    \norm{\rd\xi}_2^2 = 4n\norm{\xi}_2^2, \qquad \int\<\rd\xi,\rd\beta\> = 0. 
}
Therefore,
\eq{
    \norm{\rd\alpha}_2^2 - 4n\norm{\alpha}_2^2 = \norm{\rd\xi}_2^2 + \norm{\rd\beta}_2^2 - 4n\norm{\xi}_2^2 - 4n\norm{\beta}_2^2 =  \norm{\rd\beta}_2^2 - 4n\norm{\beta}_2^2.
}
Inserting into \eqref{eq2} gives
\eq{\label{eq3}
    E_{\kappa_*}(u) - 2n\omega_{2n+1} = \norm{\rd\beta}_2^2 - 4n\norm{\beta}_2^2 + n\Norm{\frac{*F^n}{n!}-2\alpha}_2^2.
}
Since $\beta$ is $L^2$-orthogonal to the space of Killing $1$-forms, we have
\eq{\label{eq4}
    \norm{\rd\beta}_2^2 \geq (6n+3)\norm{\beta}_2^2,
}
where $6n+3$ is the second eigenvalue of the Hodge Laplacian on the space of co-closed $1$-forms. Combining \eqref{eq3} and \eqref{eq4} gives
\eq{\label{eq5}
    \norm{\rd\beta}_2^2 + \norm{\beta}_2^2 + \Norm{\frac{*F^n}{n!}-2\alpha}_2^2 \lesssim E_{\kappa_*}(u) - 2n\omega_{2n+1}.
}

We now claim that $\xi$, $\rd\xi$, and $(\rd\xi)^n$ are all uniformly bounded by a constant depending only on $n$. In fact, the orthogonal decomposition $\alpha=\xi+\beta$ implies that
\eq{
    \norm{\xi}_2^2 \leq \norm{\alpha}_2^2 \leq \frac{1}{4n}\norm{F}_2^2 \leq \frac{1}{4n}E_{\kappa_*}(u) \leq \frac{1}{4n}(1+2n\omega_{2n+1}),
}
where we used \eqref{eq1.2} in the last step. Then by \eqref{eq2.1}, we have
\eq{
    \norm{\rd\xi}_2^2 \leq 1+2n\omega_{2n+1}.
}
Finally, since $\norm{\rd(\cdot)}_\infty$ is a norm on the space of Killing $1$-forms, which is of finite dimension, we have
\eq{
    \norm{\rd\xi}_\infty \lesssim \norm{\xi}_2.
}
By Cauchy--Schwarz,
\eq{
    \norm{(\rd\xi)^n}_2 \lesssim \norm{\rd\xi}_\infty^{n-1}\norm{\rd\xi}_2 \lesssim \norm{\xi}_2^n.
}
Hence the claim follows.

If $\abs{\rd\beta}\leq1$, then expansion gives
\eq{
    \Abs{F^n-(\rd\xi)^n}^2 = \Abs{ \sum_{j=1}^n \binom{n}{j}(\rd\xi)^{n-j}\w(\rd\beta)^j }^2 \lesssim \abs{\rd\beta}^2,
}
where we used $\abs{\rd\beta}^j\leq\abs{\rd\beta}$ in the last step.

If $\abs{\rd\beta}>1$, then
\eq{
    *\big(F^n-(\rd\xi)^n\big) = n!\Big( \frac{*F^n}{n!}-2\alpha \Big) + 2n!\beta + 2n!\xi - *(\rd\xi)^n,
}
where
\eq{
    2n!\xi - *(\rd\xi)^n \lesssim 1 < \abs{\rd\beta}.
}
Hence
\eq{
    \abs{F^n-(\rd\xi)^n}^2 \lesssim \Abs{\frac{*F^n}{n!}-2\alpha}^2 + \abs{\beta}^2 + \abs{\rd\beta}^2.
}

In both cases, combining with \eqref{eq5} gives
\eq{\label{eq5.1}
    \norm{\rd\beta}_2^2 + \norm{F^n-(\rd\xi)^n}_2^2 \lesssim E_{\kappa_*}(u) - 2n\omega_{2n+1}.
}
Next, note that $*(\rd\xi)^n$ is still a Killing $1$-form, hence
\eq{\label{eq6}
    \int \xi\w\rd\beta\w(\rd\xi)^{n-1} = \int \beta\w(\rd\xi)^n = \int\<\beta,*(\rd\xi)^n\> = 0,
}
where we used Stokes' theorem in the first step. Moreover, if $\abs{\rd\beta}\leq1$, then expansion gives
\eq{
    \Abs{ F^n - (\rd\xi)^n - n\rd\beta\w(\rd\xi)^{n-1} } = \Abs{ \sum_{j=2}^n \binom{n}{j}(\rd\xi)^{n-j}\w(\rd\beta)^j } \lesssim \abs{\rd\beta}^2.
}
If $\abs{\rd\beta}>1$, then
\eq{
    \Abs{ F^n - (\rd\xi)^n - n\rd\beta\w(\rd\xi)^{n-1} } \leq \Abs{ F^n - (\rd\xi)^n } + C(n)\abs{\rd\beta} \lesssim \Abs{ F^n - (\rd\xi)^n }^2 + \abs{\rd\beta}^2.
}
In both cases, combining with \eqref{eq6} gives
\eq{\label{eq7}
    \Abs{ \int \xi\w\big(F^n-(\rd\xi)^n\big) } &= \Abs{ \int \xi\w\big(F^n - (\rd\xi)^n - n\rd\beta\w(\rd\xi)^{n-1}\big) } \\
    &\lesssim \norm{\rd\beta}_2^2 + \norm{F^n-(\rd\xi)^n}_2^2 \lesssim E_{\kappa_*}(u) - 2n\omega_{2n+1},
}
where we used \eqref{eq5.1} in the last step. It follows that
\eq{\label{eq8}
    \frac{n!}{2}\omega_{2n+1} = \int \alpha\w F^n = \int \xi\w(\rd\xi)^n + \int \xi\w\big(F^n-(\rd\xi)^n\big) + \int \beta\w F^n.
}
Here
\eq{\label{eq9}
    \Abs{\int \beta\w F^n} = \Abs{\int \beta\w \big(F^n - (\rd\xi)^n\big)} \lesssim \norm{\beta}_2\cdot\Norm{F^n - (\rd\xi)^n}_2 \lesssim \norm{\beta}_2^2 + \Norm{F^n - (\rd\xi)^n}_2^2 \lesssim E_{\kappa_*}(u) - 2n\omega_{2n+1},
}
where we used \eqref{eq6} in the first step and \eqref{eq5.1} in the last step. Combining \eqref{eq7}, \eqref{eq8}, and \eqref{eq9} gives
\eq{\label{eq10}
    \Abs{ \frac{2}{n!\omega_{2n+1}}\int\xi\w(\rd\xi)^n - 1 } \lesssim E_{\kappa_*}(u) - 2n\omega_{2n+1}.
}
Recall that $\xi$, $\rd\xi$, and $(\rd\xi)^n$ are all uniformly bounded, hence \eqref{eq5} and \eqref{eq5.1} imply that
\eq{\label{eq10.1}
    \Norm{ \frac{*(\rd\xi)^n}{n!}-2\xi }_2^2 = \Norm{ \Big( \frac{*F^n}{n!}-2\alpha \Big) - \frac{*(F^n-(\rd\xi)^n)}{n!} + 2\beta }_2^2 \lesssim E_{\kappa_*}(u) - 2n\omega_{2n+1}.
}

Since the right-hand side of \eqref{eq1.1} is bounded by $E_{\kappa_*}(u)$, in order to prove \eqref{eq1.1}, we only need to consider the case that $E_{\kappa_*}(u) - 2n\omega_{2n+1}$ is sufficiently small. Then \eqref{eq10} implies that
\eq{
    \int \xi\w(\rd\xi)^n > 0,
}
which means $\xi$ has no zero singular values. Therefore, we may write in suitable coordinates that
\eq{\label{eq10.1.1}
    \xi = \frac{1}{2}\sum_{j=1}^{n+1} r_j(x_j\rd y_j-y_j\rd x_j), \qquad r_j>0,
}
and hence
\eq{\label{eq10.1.2}
    \rd\xi = \sum_{j=1}^{n+1} r_j\,\rd x_j\w\rd y_j.
}
It follows that
\eq{\label{eq10.1.3}
    \frac{*(\rd\xi)^n}{n!} = \sum_j \Big( \prod_{k\not=j}r_k \Big) (x_j\rd y_j-y_j\rd x_j).
}
Hence
\eq{\label{eq10.2}
    \Norm{ \frac{*(\rd\xi)^n}{n!}-2\xi }_2^2 &= \Norm{ \sum_j \Big( \prod_{k\not=j}r_k - r_j \Big) (x_j\rd y_j-y_j\rd x_j) }_2^2 \\
    &= \sum_j \Big( \prod_{k\not=j}r_k - r_j \Big)^2 \Big( \norm{x_j}_2^2 + \norm{y_j}_2^2 \Big) \\
    &= \frac{\omega_{2n+1}}{n+1}\sum_j \Big( \prod_{k\not=j}r_k - r_j \Big)^2,
}
where we used the fact that $\norm{x_j}_2^2=\norm{y_j}_2^2=\frac{\omega_{2n+1}}{2n+2}$. Moreover, we have
\eq{\label{eq11}
    \int \xi\w(\rd\xi)^n = \int \<\xi,*(\rd\xi)^n\> = \frac{n!}{2}\Big( \prod_j r_j \Big) \sum_j \big(\norm{x_j}_2^2+\norm{y_j}_2^2\big) = \frac{n!}{2}\omega_{2n+1}\Big( \prod_j r_j \Big).
}

Note that for each $j$ we have
\eq{
    r_j^2-1 = -r_j\Big( \prod_{k\not=j}r_k-r_j \Big) + \Big( \prod_jr_j-1 \Big),
}
and
\eq{
    (r_j-1)^2 &= \Big( \frac{r_j^2-1}{r_j+1} \Big)^2 = \Big\{ -\frac{r_j}{r_j+1}\Big( \prod_{k\not=j}r_k-r_j \Big) + \frac{1}{r_j+1}\Big( \prod_jr_j-1 \Big) \Big\}^2 \\
    &\leq 2\Big( \prod_{k\not=j}r_k-r_j \Big)^2 + 2\Big( \prod_jr_j-1 \Big)^2.
}
It follows that
\eq{\label{eq12}
    \sum_j (r_j-1)^2 \leq 2\sum_j\Big( \prod_{k\not=j}r_k-r_j \Big)^2 + 2(n+1)\Big( \prod_jr_j-1 \Big)^2 \lesssim E_{\kappa_*}(u) - 2n\omega_{2n+1},
}
where we used \eqref{eq10.1} and \eqref{eq10.2} to estimate the first term, and used \eqref{eq10} and \eqref{eq11} for the second.

Expanding at $(r_1,\cdots,r_{n+1})=(1,\cdots,1)$ gives
\eq{
    \prod_j r_j - 1 = \sum_j(r_j-1) + O_n\Big( \sum_j(r_j-1)^2 \Big),
}
hence
\eq{\label{eq13}
    \sum_j (r_j^2-1) = 2\sum_j(r_j-1) + \sum_j(r_j-1)^2 = 2\Big( \prod_j r_j - 1 \Big) + O_n\Big( \sum_j(r_j-1)^2 \Big).
}
Combining \eqref{eq10}, \eqref{eq11}, \eqref{eq12} and \eqref{eq13} gives
\eq{\label{eq14}
    \abs{\sum_j (r_j^2-1)} \lesssim E_{\kappa_*}(u) - 2n\omega_{2n+1}.
}

Finally, \eqref{eq10.1.1} implies that
\eq{
    \norm{\xi}_2^2 = \frac{\omega_{2n+1}}{4(n+1)}\sum_jr_j^2,
}
and hence
\eq{
    \norm{\rd\xi}_2^2 = 4n\norm{\xi}_2^2 = \frac{n\omega_{2n+1}}{n+1}\sum_jr_j^2 = n\omega_{2n+1} + \frac{n\omega_{2n+1}}{n+1}\sum_j (r_j^2-1).
}
Using \eqref{eq5} and \eqref{eq14} we have
\eq{\label{eq15}
    \norm{F}_2^2 = \norm{\rd\xi}_2^2 + \norm{\rd\beta}_2^2 = n\omega_{2n+1} + O_n\big(E_{\kappa_*}(u) - 2n\omega_{2n+1}\big),
}
and hence
\eq{\label{eq16}
    \kappa_*\norm{F^n}_2^2 = E_{\kappa_*}(u) - \norm{F}_2^2 &= E_{\kappa_*}(u) - n\omega_{2n+1} + O_n\big(E_{\kappa_*}(u) - 2n\omega_{2n+1}\big) \\
    &= n\omega_{2n+1} + O_n\big(E_{\kappa_*}(u) - 2n\omega_{2n+1}\big).
}
Combining \eqref{eq15} and \eqref{eq16} gives
\eq{
    \Abs{ \norm{F}_2^2 - \kappa_*\norm{F^n}_2^2 } \lesssim E_{\kappa_*}(u) - 2n\omega_{2n+1},
}
which completes the proof of \eqref{eq1.1}.
\end{proof}

We now prove Theorem~\ref{thm4.1}.

\begin{proof}[Proof of Theorem~\ref{thm4.1}]
Note that
\eq{\label{eq17}
    E_\kappa(u) - E_\kappa(\pi) = \frac{\kappa+\kappa_*}{2\kappa_*}(E_{\kappa_*}(u)-2n\omega_{2n+1}) + \frac{\kappa_*-\kappa}{2\kappa_*}\Big( \int\abs{F}^2 - \kappa_*\int\abs{F^n}^2 \Big).
}
Combining \eqref{eq1.1} and \eqref{eq17} gives
\eq{\label{eq18}
    E_\kappa(u) - E_\kappa(\pi) \geq \Big( \frac{\kappa+\kappa_*}{2\kappa_*} - \frac{\abs{\kappa_*-\kappa}}{2\kappa_*}\frac{1}{2\epsilon_n} \Big)(E_{\kappa_*}(u)-2n\omega_{2n+1}).
}
By \eqref{eq1.0}, we see that the condition \eqref{eq1} is equivalent to
\eq{
    \frac{1-\epsilon_n}{1+\epsilon_n}\kappa_* \leq \kappa \leq \frac{1+\epsilon_n}{1-\epsilon_n}\kappa_* \quad\iff\quad \abs{\kappa-\kappa_*}\leq\epsilon_n(\kappa+\kappa_*).
}
Then \eqref{eq18} implies that
\eq{\label{eq19}
    E_\kappa(u) - E_\kappa(\pi) \geq \frac{\kappa+\kappa_*}{4\kappa_*}(E_{\kappa_*}(u)-2n\omega_{2n+1}) \geq 0.
}
Hence the Hopf map $\pi$ minimizes $\mathcal{E}_{{\rm RTY}}^R$ whenever \eqref{eq1} holds.

Assume equality holds. Then \eqref{eq19} implies that $E_{\kappa_*}(u)=2n\omega_{2n+1}$. By \eqref{eq2}, we see that
\eq{
    \norm{\rd\alpha}_2^2 = 4n\norm{\alpha}_2^2, \qquad \frac{*(\rd\alpha)^n}{n!}=2\alpha,
}
which means $\alpha$ is a Killing $1$-form, hence $\alpha=\xi$ and $\beta=0$. Recall \eqref{eq10.1.1} that in suitable coordinates
\eq{\label{eq20}
    \alpha = \frac{1}{2}\sum_{j=1}^{n+1} r_j(x_j\rd y_j-y_j\rd x_j).
}
Comparing \eqref{eq11} with $\int \alpha\w(\rd\alpha)^n=\frac{n!}{2}\omega_{2n+1}$ gives
\eq{\label{eq21}
    \prod_jr_j=1.
}
Comparing \eqref{eq10.1.3} with \eqref{eq20} gives
\eq{\label{eq22}
    \prod_{k\not=j}r_k = r_j, \qquad\forall\,j.
}
By \eqref{eq21} and \eqref{eq22} we see that each $r_j=1$. Therefore, up to a $\mathcal{R}\in{\rm SO}(2n+2)$, we have
 \eq{
     \alpha = \mathcal{R}^*\eta, \qquad \eta = \frac{1}{2}\sum_{j=1}^{n+1} (x_j\rd y_j-y_j\rd x_j), \qquad \rd\eta=\pi^*\omega.
 }
 Hence
\eq{
    u^*\omega = \rd\alpha = \mathcal{R}^*\rd\eta = \mathcal{R}^*\pi^*\omega = (\pi\circ\mathcal{R})^*\omega.
}
For $v\coloneqq u\circ\mathcal{R}^{-1}$, we have
\eq{
    v^*\omega = \pi^*\omega.
}
Since $\rd\pi$ has rank $2n$ at every point, so is $\rd v$, and
\eq{
    \ker(\rd v) = \ker(v^*\omega) = \ker(\rd\pi).
}
That is, $v$ is invariant along each connected Hopf $\S^1$-fibre. Hence it descends to
\eq{
    h:\CP^n\to\CP^n, \qquad v=h\circ\pi.
}
Since
\eq{
    \pi^*(h^*\omega-\omega) = v^*\omega - \pi^*\omega = 0,
}
we see that $h^*\omega=\omega$. Moreover, we have
\eq{
    \deg(h)\int_{\CP^n}\omega^n = \int_{\CP^n} h^*\omega^n = \int_{\CP^n} \omega^n,
}
hence $\deg(h)=1$. Therefore, $h$ is a symplectic diffeomorphism, and $u=h\circ\pi\circ\mathcal{R}$. The reverse is trivial. Hence we complete the proof.
\end{proof}

\subsection{Proof of Theorem~\ref{thm4.2}}

\begin{proof}[Proof of Theorem~\ref{thm4.2}]

Let $\Phi_t$ be the conformal flow generated by $\nabla x$, i.e.
\eq{
    \partial_t\Phi_t = \nabla x\circ\Phi_t, \qquad \Phi_0=\rid.
}
It is clear that
\eq{
    \mathcal{L}_{\nabla x}g = -2xg, \qquad \frac{\rd}{\rd t}\Phi_t^*g = -2(x\circ\Phi_t)\Phi_t^*g.
}
The orbit $y(t)=x(\Phi_t(\cdot))$ solves
\eq{
    y'(t) = 1-y(t)^2, \qquad y(0)=x,
}
hence
\eq{
    y(t) = \frac{x\cosh t + \sinh t}{\cosh t + x\sinh t}.
}
Therefore,
\eq{
    \Phi_t^*g = \lambda_t^2g, \qquad \lambda_t = \frac{1}{\cosh t + x\sinh t}.
}
It follows that
\eq{
    \abs{\Phi_t^*(\pi^*\omega)}^2 = n\lambda_t^4, \qquad \abs{\Phi_t^*(\pi^*\omega)^n}^2 = (n!)^2\lambda_t^{4n}.
}
Hence
\eq{\label{eq23}
    E_\kappa(\pi\circ\Phi_t) = n\int \lambda_t^4 + \kappa (n!)^2\int \lambda_t^{4n}.
}

For any $p$, expansion gives
\eq{
    \lambda_t^p = 1 - pxt + \Big( \frac{p(p+1)}{2}x^2 - \frac{p}{2} \Big)t^2 + O(t^3). 
}
Since
\eq{
    \frac{1}{\omega_{2n+1}}\int x = 0, \qquad \frac{1}{\omega_{2n+1}}\int x^2 = \frac{1}{2n+2},
}
we have
\eq{
    \frac{1}{\omega_{2n+1}}\int \lambda_t^p = 1 + \frac{p(p-2n-1)}{4n+4}t^2 + O(t^4).
}
In particular, choosing $p=4$ and $p=4n$ gives
\eq{
    n\int \lambda_t^4 = n\omega_{2n+1} - \frac{n(2n-3)\omega_{2n+1}}{n+1}t^2 + O(t^4)
}
and
\eq{
    \kappa(n!)^2\int \lambda_t^{4n} = \kappa(n!)^2\omega_{2n+1} + \frac{\kappa(n!)^2n(2n-1)\omega_{2n+1}}{n+1}t^2 + O(t^4).
}
Combining with \eqref{eq1.0.1} and \eqref{eq23} gives
\eq{
    E_\kappa(\pi\circ\Phi_t) - E_\kappa(\pi) = \frac{n\omega_{2n+1}}{n+1}\big((2n-1)\kappa(n!)^2-(2n-3)\big)t^2 + O(t^4).
}
Therefore, $\pi$ is not locally stable whenever
\eq{
    \kappa < \frac{2n-3}{(2n-1)(n!)^2},
}
which is exactly \eqref{eq1.-2} in view of \eqref{eq1.0}. Since the generating vector fields, i.e. the conformal Killing vector fields that are the gradient fields of coordinate functions, have dimension $2n+2$, we have ${\rm ind}(\pi)\geq2n+2$.

We now look at the other side. Choose a horizontal vector field $X$ on $\S^{2n+1}$ such that $\rdiv X=0$ and $[X,\xi]=J_HX$, where $\xi$ is the standard Reeb field and $J_H$ is the horizontal complex structure. In suitable complex coordinates, we consider
\eq{
    X \coloneqq (-iz_3\bar{z_2},iz_3\bar{z_1},0,\cdots,0).
}
Set
\eq{
    \sigma \coloneqq X\ip(\pi^*\omega) = \rRe[\bar{z_3}(z_2\rd z_1-z_1\rd z_2)].
}
Direct computation checks that
\eq{\label{eq23.1}
    \norm{\rd\sigma}_2^2 = (6n+3)\norm{\sigma}_2^2.
}
Let $\phi_t$ be the flow generated by $X$, then $\phi_t$ preserves the volume since ${\rm div}X=0$. The first variation is
\eq{
    \frac{\rd}{\rd t}\Big|_{t=0} \phi_t^*(\pi^*\omega) = \rd\sigma, \qquad \frac{\rd}{\rd t}\Big|_{t=0} \phi_t^*\xi = J_HX, 
}
and the second variation is
\eq{
    \frac{\rd^2}{\rd t^2}\Big|_{t=0} \phi_t^*(\pi^*\omega) = \rd(X\ip\rd\sigma), \qquad \frac{\rd^2}{\rd t^2}\Big|_{t=0} \phi_t^*\xi = \mathcal{L}_X(J_HX) = [X,J_HX], 
}
By $\rd^*(\pi^*\omega)=2n\xi^\flat$, it is clear that $\pi$ is a critical point of $E_\kappa$ for each $\kappa>0$. Moreover, we have
\eq{\label{eq24}
    \frac{1}{2}\frac{\rd^2}{\rd t^2}\Big|_{t=0} \norm{\phi_t^*(\pi^*\omega)}_2^2 = \norm{\rd\sigma}_2^2 + \int \<2n\xi^\flat,X\ip\rd\sigma\> = \norm{\rd\sigma}_2^2 - 2n\norm{\sigma}_2^2 = (4n+3)\norm{\sigma}_2^2.
}
Using $\phi_t^*(\pi^*\omega)^n=\phi_t^*(n!\xi\ip\rdV)=n!(\phi_t^*\xi)\ip\rdV$ we have
\eq{\label{eq25}
    \frac{1}{2}\frac{\rd^2}{\rd t^2}\Big|_{t=0} \norm{\phi_t^*(\pi^*\omega)^n}_2^2 &= (n!)^2\int\Big( \abs{J_HX}^2 + \<\xi,[X,J_HX]\> \Big) = (n!)^2\int\Big( \abs{J_HX}^2 - \rd\xi^\flat(X,J_HX) \Big) \\
    &= (n!)^2\int\Big( \abs{J_HX}^2 - 2\abs{X}^2 \Big) = -(n!)^2\norm{\sigma}_2^2.
}
Combining \eqref{eq24} and \eqref{eq25} gives
\eq{
    E_\kappa(\pi\circ\phi_t) - E_\kappa(\pi) = \big(4n+3-\kappa(n!)^2\big)\norm{\sigma}_2^2\,t^2 + O(t^3).
}
Therefore, $\pi$ is not locally stable whenever
\eq{
    \kappa > \frac{4n+3}{(n!)^2},
}
which is exactly \eqref{eq1.-1} in view of \eqref{eq1.0}. To count the dimension of the space of generating vector fields, it is convenient to assume that $X$ is quadratic, i.e.
\eq{
    X_i(z) = \sum_{j,k=1}^{n+1} c_{ijk}\bar{z_j}z_k,
}
where the form reflects that $[X,\xi]=J_HX$ (which helps to decrease the second variation). The conditions imposed on $X$ then become
\eq{
    c_{ijk}=-c_{jik}, \qquad \sum_j c_{ijj}=0,
}
where the first reflects that $X$ is horizontal, and the second reflects that ${\rm div}X=0$. The corresponding $\sigma\coloneqq X\ip(\pi^*\omega)$ satisfies \eqref{eq23.1}. It is then easy to see that the complex dimension of these $X$'s is
\eq{
    (n+1)\binom{n+1}{2}-(n+1) = \frac{1}{2}(n-1)(n+1)(n+2).
}
Therefore ${\rm ind}(\pi)\geq(n-1)(n+1)(n+2)$.
Hence we complete the proof.
\end{proof}

\printbibliography

\end{document}